\documentclass[11pt]{amsart}
\usepackage{graphicx}
\usepackage{amssymb}
\usepackage{amsmath,amssymb,amsthm,epsfig}

\newtheorem{theorem}{Theorem}[section]
\newtheorem{proposition}[theorem]{Proposition}
\newtheorem{lemma}[theorem]{Lemma}

\newtheorem{definition}[theorem]{Definition}

\title{Bounding right-arm rotation distances}
\author{Sean Cleary}
\address{Department of Mathematics,
The City College of New York \& The CUNY Graduate Center, New York, NY 10031} 
\email{cleary@sci.ccny.cuny.edu}
\thanks{The first author acknowledges support from PSC-CUNY grant \#66490,
NSF grant DMS-0305545
and the hospitality of the Centre de Recerca Matem\`atica.}
\author{Jennifer Taback}
\address{Department of Mathematics,
Bowdoin College, Brunswick, ME 04011} \email{jtaback@bowdoin.edu}
\thanks{The second author acknowledges support from
NSF grants DMS-0305545 and DMS-0437481 and the hospitality of the Centre de Recerca Matem\`atica.
}

\begin{document}
\begin{abstract}
Rotation distance measures the difference in shape between
 binary trees of the same size by counting the minimum
number of rotations needed to transform one tree to the
other.  We describe several types of rotation distance where restrictions
are put on the locations where rotations are permitted, and
provide upper bounds  on distances between trees with a fixed number
of nodes with respect to several families of these restrictions.  These bounds are
sharp in a certain asymptotic sense and are obtained by relating each restricted rotation distance
to the word length of
elements of Thompson's group $F$ with respect to different generating
sets, including both finite and infinite generating sets.
\end{abstract}
\maketitle

\section{Introduction}

Rotation  distance quantifies the difference in shape between two
rooted binary trees of the same size by counting the minimum
number of elementary changes needed to transform one tree to the
other.  Search algorithms are most efficient when
searching balanced trees, which have few levels relative to the
number of nodes in the tree.  Thus one is often interested in
calculating, or at least bounding, the number of these changes
necessary to alter a given tree into another with a more desirable shape,
such as a balanced tree.

If  we allow these elementary changes, called rotations, to take place at
any node, we obtain ordinary rotation distance.  This was
analyzed by Sleator, Tarjan and Thurston \cite{stt}, who proved an upper bound of $2n-6$ rotations needed to transform one rooted
binary tree with $n$ nodes into any other,  for $n \geq 11$.  Furthermore,
they showed that the $2n-6$ bound is achieved for 
all sufficiently large $n$ and thus is the best possible upper bound. No efficient
algorithm is known to compute rotation distance exactly, though
there are polynomial-time algorithms of Pallo \cite{pallo} and
Rogers \cite{rogers} which estimate rotation distance efficiently.

Here we expand on the study of restricted rotation distance begun in \cite{rotipl}
and \cite{rotbound}.  Restricted rotation distance allows rotations only at the root node and the right child of the root node.  Restricted rotation distance is related to the word length of elements of
Thompson's group $F$ with respect to its standard finite generating set.  This is illustrated in \cite{rotipl, rotbound} and involves the interpretation of elements of $F$ as pairs of finite binary rooted trees and
Fordham's method  \cite{blake:gd} for computing the word length of an element of $F$ with respect to
that standard finite generating set directly from such trees.  These methods not only
give an effective algorithm to compute restricted rotation distance, but they also give
an effective algorithm to find the appropriate rotations which realize this distance.

\begin{figure}
\includegraphics[width=2.5in]{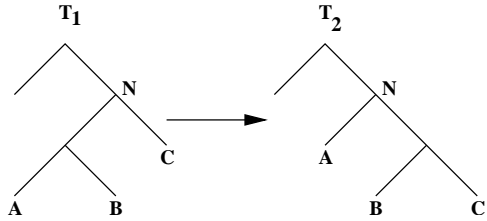}
\caption{Right rotation at node $N$ transforms $T_1$ to $T_2$;
similarly, left rotation at node $N$ transforms $T_2$ to $T_1$.
The labels $A$, $B$ and $C$ represent (possibly empty) subtrees of $T_1$ and $T_2$.
\label{fig:rightfig1}}
\end{figure}

Right and left rotations at a node $N$ of a rooted binary tree $T$  are
defined to be the permutations of the subtrees of $T$  described in
Figure \ref{fig:rightfig1}.  Right rotation at a node $N$
transforms the original tree $T_1$, given on the left side of
Figure \ref{fig:rightfig1}, to the tree $T_2$ on the right side of
Figure \ref{fig:rightfig1}. Left rotation at a node is the inverse
operation. In all that follows, $T_1$ and $T_2$ denote trees with
the same number of nodes.

In this paper, we discuss generalizations and variations of restricted
rotation distance, in which rotations are again only allowed at
specified nodes of the tree.  
We relate these  distances to 
distinct word metrics on Thompson's group $F$.  We use this
interpretation to exhibit linear bounds on the number of allowable elementary
rotations needed to transform one tree with $n$
nodes into another, and show that these bounds are asymptotically sharp in the sense that the
coefficients of the linear terms of the bounds are the best possible.  These
 alternate definitions all
allow rotations at the root node and at nodes connected to the root node by a path
consisting entirely of right edges; that is, nodes that lie on the 
{\em right side} or {\em right arm} of the tree.  The root node is considered to lie on the right side of the tree. 

One complication that arises is that while the original restricted rotation distance is
always defined between any two trees with the same number of nodes, this is no longer
necessarily the case when we allow rotations at  other 
collections of nodes along the right side of the tree.
Some transformations between trees cannot be accomplished with
a specified set of rotations  without the nonstandard technique of
adding additional nodes to the trees.  Such a transformation
is not permitted when computing rotation distances of any type.
Below, we describe how to determine when such transformations are possible with a prescribed set of permitted locations at which to rotate. When this restricted right-arm rotation distance is defined, we provide an upper bound on its magnitude.

The sharp upper bound on the restricted rotation distance between two trees, each with $n$ nodes, obtained in \cite{rotbound} is $4n-8$, for $n \geq 3$.
Below, we consider allowing additional rotations along the right 
side of the tree and note that allowing rotations at any finite collection of nodes
on the right side of the tree does not change the multiplicative
constant of 4 in the upper bound.  It is only when we allow an
infinite set of rotations along the right arm of the tree that we
obtain the multiplicative constant of 2 in the upper bound, analogous to ordinary rotation distance.
These rotation distances and bounds, which hold for sufficiently large $n$, are summarized in Table \ref{tab:summ}, where $n$ is the number of nodes in each tree.

 \begin{table}[h]\label{summarytable}
 \begin{tabular}{|l|l|c|c|}
 \hline
 {\bf Type of distance} &  {\bf Rotations allowed at} & {\bf Symbol} & {\bf Upper Bound}\\
 \hline
 Rotation distance & all nodes & $d_R$ & $2n-6$\\
\hline
 Restricted rotation & root node and right & $d_{RR}$ & $4n-8$ \\
 distance  & child of the root node & & \\ \hline
 Restricted right arm  & root node and  a finite&  &
 $4n-C$\\
 rotation distance & collection $\mathcal{S}$ of nodes on & $d_{RRA}^\mathcal{S}$ &  some $C$ \\
 & the right side of the tree & &  \\
\hline
 Right arm  & all nodes on the right  & $d_{RA}$ & $2n-2$\\
 rotation distance & side of the tree & & \\ \hline
 Restricted spinal  & root node and a finite& &
 $4n-C$\\
 rotation distance & collection $\mathcal{S}$ of nodes on  & $d_{RS}^\mathcal{S}$ & some $C$ \\
 & both sides of the tree & &  \\
 \hline 
\end{tabular}
\vspace{.1in}
\caption{Summary of rotation distances between trees with $n$ nodes, with upper bounds for all sufficiently large $n$. \label{tab:summ}}
\end{table}

Culik and Wood \cite{cw}, in the course of studying ordinary rotation distance,
showed that the rotation distance is never more than $2n-2$, and in fact
use only rotations on the right arm to show this bound.
Since
the ordinary rotation distance between two such trees can only be as
much as $2n-6$ for $n \geq 11$, it is remarkable that restricting rotations to the
right side of the tree adds only four rotations to the upper
bound. 

Pallo explicitly studied  right-arm rotation distance in
\cite{palloright}, allowing rotations at all nodes along the right side of the tree.  He
described an algorithm for computing right-arm rotation distance which we show below
is equivalent to finding the word length in Thompson's group $F$ with respect to
the standard infinite generating set. 

The trees we consider are composed of edges and vertices.  The vertices fall into two types: those of valence one and those of higher valence.  The vertices of valence one are called {\em exterior nodes} or {\em leaves} or {\em exposed leaves}.  The vertices of higher valence are called {\em interior nodes}.  We shall use the terminology {\em node} to refer to a vertex which is an interior node, and {\em leaf} to refer to a vertex which is an exterior node.

A {\em caret} in a tree is composed of a node together with two downward directed edges.  We will only consider finite, rooted binary trees with $n$ carets, equivalently, with $n$ nodes.  Such trees are called {\em extended binary trees} in Knuth \cite{knuth3} or {\em 0-2 trees}.  The nodes and carets in a tree have a natural infix ordering.  The exposed leaves in a tree are numbered from left to right,
beginning with zero.  A tree with $n$ carets yields $n+1$ exposed leaves.  
A caret with two exposed leaves is called an {\em exposed
caret}, its leaves are termed {\em siblings} and those leaves
are said to form a {\em sibling pair}.

A caret $N$ which is attached to the right (respectively  left) edge of a caret $M$ is called the {\em right (respectively left) child} of $M$.  A caret which has one edge on the left side of the tree is called a {\em left caret}.  A caret which has one edge on the right side of the tree and is not the root caret is called a {\em right caret}. Similarly, we have left and right nodes.
Carets which are neither right nor left  are called {\em interior carets}.  The union of left and right carets in a tree is called the {\em spine} of the tree.  A tree consisting of only the root caret and $n-1$ right
carets is called the {\em all-right tree} with $n$ carets.  An {\em ancestor} of a caret (resp. node) is any caret (resp. node) which lies along the shortest path between it and the root caret (resp. root node).

The connection between Thompson's group $F$ and
restricted rotation distance is described below.  Thompson's group $F$
is studied combinatorially in two ways: via a finite presentation and an
infinite presentation.  Computing restricted rotation
distance between two trees is related to computing the word length of the element of $F$ described by those trees with respect to the standard finite generating set for the group $F$.  Analogously,
right-arm rotation distance corresponds to computing the word length of the element with respect to the word metric induced by the standard infinite generating set for $F$.
Restricted right-arm rotation distances and restricted spinal
rotation distances, defined below, relate to the
word metric on $F$ with respect to other finite generating sets.

\section{Thompson's Group $F$}

The connection between Thompson's group $F$ and rotations at nodes
of trees is described in \cite{rotipl} and \cite{rotbound}, using
the work of Fordham \cite{blake:gd}. Here, we briefly describe this
connection, and refer the reader to Cannon, Floyd and Parry
\cite{cfp} for a survey of the properties of Thompson's group $F$, and the 
further connections between elements of $F$ and pairs of
binary rooted trees.

\subsection{The infinite presentation of Thompson's group $F$}

Thompson's group $F$ has a presentation with an infinite number of
generators and relations: $$ {\mathcal P} = \left<x_0,x_1,\ldots |
x_i^{-1} x_n x_i=x_{n+1},\forall i < n \right>.$$  In this
presentation, there are normal forms for elements given by
$$x_{i_1}^{r_1} x_{i_2}^{r_2}\cdots x_{i_k}^{r_k} x_{j_l}^{-s_l}
\cdots x_{j_2}^{-s_2} x_{j_1}^{-s_1} $$ with $r_i, s_i >0$, where the indices satisfy $0 \leq i_1<i_2 < \cdots < i_k$ and $0 \leq j_1<j_2 < \cdots < j_l$. This normal
form is unique for a given element if we further require the
reduction condition that when both $x_i$ and $x_i^{-1}$ occur, so
does $x_{i+1}$ or  $x_{i+1}^{-1}$, as discussed by Brown and
Geoghegan \cite{bg:thomp}. The relators provide a quick and
efficient method for rewriting words into normal form, and form a
complete rewriting system, as described by Brown \cite{brownfws}.
There is a natural shift homomorphism $\phi:F \rightarrow F$
where $\phi(x_i)=x_{i+1}$ which respects the relators, and
the reduction from normal form  to unique normal form is accomplished
with a sequence of operations replacing words of the form
$u x_i \phi(v) x_i^{-1} w$ with $u v w$, where $\phi(v)$ is a subword
which contains only generators  of index ${i+2}$ and higher.

We note that $F$ can be generated by just $x_0$ and
$x_1$ in the above presentation; the relators  show that $x_0$
conjugates $x_1$ to $x_2$. Similarly,  all higher-index generators are
conjugates of $x_1$ by higher powers of $x_0$, as $x_n= x_0^{-(n-1)} x_1 x_0^{n-1}$.  This leads to a
finite presentation for $F$ with generating set $\{x_0,x_1\}$.  In
fact, $x_0$ and any higher index generator are sufficient to
generate the group.  Any two generators ${x_i, x_j}$ with $i \neq j$ will
generate an subgroup of $F$ which is isomorphic to the entire group
but which is the entire group only when one of $i$ or $j$ is $0$.

We begin by proving that in the word metric arising from this infinite generating set, the normal form expressions are geodesic representatives for elements of $F$.

\begin{lemma}\label{lemma:geodesic}
Let $w$ be an element of $F$, and $\alpha$ a word in the infinite generating set which is the unique normal form for $w$, as described above.  Then $\alpha$ is a geodesic representative for $w$ in the word metric
arising from the infinite generating set $\{x_i\}$ of $F$.
\end{lemma}

\begin{proof}
Suppose that $\alpha = x_{i_1}^{r_1} x_{i_2}^{r_2}\cdots x_{i_k}^{r_k}
x_{j_l}^{-s_l} \cdots x_{j_2}^{-s_2} x_{j_1}^{-s_1} $ was not a
geodesic representative for $w$ in this word metric. Then there is a shorter
expression $\beta$, not necessarily in normal form, representing $w$ in this infinite generating set.  It is clear from the relations of ${\mathcal P}$ that the conversion of $\beta$ into unique normal form can only  preserve or
decrease the length of $\beta$.  Thus, after converting $\beta$ into normal form we have obtained a second expression for $w$ in unique normal form shorter than the initial unique normal form for $w$ given by $\alpha$,  a contradiction.
 \end{proof}

\subsection{Tree pair diagrams for elements of Thompson's group $F$}

The group $F$ has a geometric description in terms of 
equivalence classes of tree pair
diagrams. A {\em tree pair diagram} is a pair of finite rooted binary
trees with the same number of nodes (or carets), or equivalently with the same number of leaves.
 We write $w = (T_1,T_2)$ to
denote the two trees comprising a pair representing $w$.  The
equivalence between the geometric and algebraic
interpretations of $F$ is described in \cite{cfp}, and examples of
this equivalence and its connection with rotations are given in \cite{ct2}.

\begin{figure}
\includegraphics[width=2.5in]{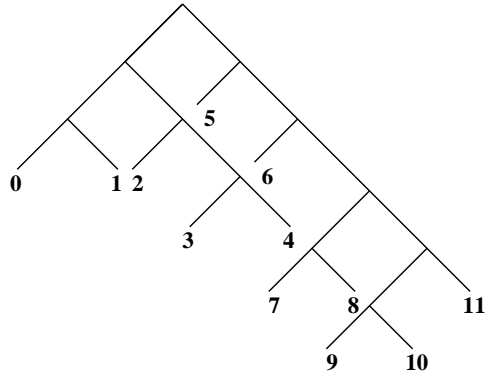} \\
\caption{A tree whose leaves are numbered from left to right.  The leaf exponents of the leaves, according to increasing leaf number, are 2,0,1,1,0,0,0,1,0,1,0, and 0.}
\label{fig:expexample}
\end{figure}

Given two trees $T_1$ and $T_2$ with the same number of carets, the word in normal form associated to $w = (T_1,T_2)$ is found as follows.  The leaves of each tree are numbered from left to right, beginning with zero.
The {\em leaf exponent} of a leaf numbered $k$ is the integral
 length of the longest path starting at leaf $k$ consisting entirely of left
 edges which does not touch the right side of the tree.  The tree pair
 diagram $(T_1,T_2)$ has an associated normal form 
 $x_0^{f_0} x_1^{f_1} \cdots x_n^{f_n} x_n^{-e_n} \cdots x_1^{-e_1} x_0^{-e_0}$ where
 $e_i$ is the leaf exponent of leaf $i$ in tree $T_1$ and $f_i$ is the leaf exponent of leaf
 $i$ in $T_2$.  
 An example of a tree with leaf exponents computed
 is given in Figure \ref{fig:expexample}.

An element of $F$ is represented uniquely by a tree pair diagram satisfying the following reduction condition.  A tree pair diagram $(T_1,T_2)$ is {\em unreduced} if both $T_1$  and $T_2$ contain a caret with two exposed leaves numbered $i$ and
$i+1$. A tree pair diagram which is not unreduced is {\em reduced}.
Geometrically, any tree pair diagram has a unique reduced form that is obtained by
successively deleting exposed carets with
identical leaf numbers from both trees, renumbering the leaves, and repeating this process until no further such reductions are possible.  Elements of $F$ are
 equivalence classes of tree pair diagrams, where the equivalence relation
 is that two tree pairs are equivalent if they have a common reduced form.

This tree pair reduction condition corresponds exactly to the
combinatorial reduction condition given above to ensure uniqueness
for words in normal form in the infinite presentation of $F$. 
That is, if leaves $i$ and $i+1$ form a sibling pair in both $T_1$ and $T_2$,
then in both cases, the leaf exponents of leaf $i$ will be non-zero  in both trees and those for
leaf $i+1$ will be zero, as leaf $i+1$ is a right leaf in both trees.  So the corresponding
normal form will contain both $x_i$ and $x_i^{-1}$ but not $x_{i+1}^{\pm1}$, meaning that the normal form can be reduced.

To perform the group operation on the level of tree pair diagrams, it may be necessary to use unreduced representatives of elements.  Namely, to multiply $(T_1,T_2)$ and $(S_1,S_2)$, we create unreduced representatives $(T'_1,T'_2)$ and $(S'_1,S'_2)$ in which $T'_2 = S'_1$, and write the product as the (possibly unreduced) element $(T'_1,S'_2)$.  See \cite{cfp} for examples of group multiplication using tree pair diagrams for elements of $F$.

\begin{figure}
\includegraphics[width=5in]{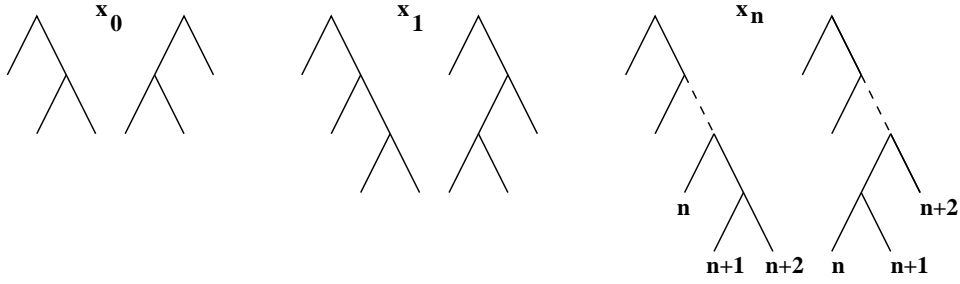} \\
\caption{The tree pair diagrams corresponding to the generators $x_0$, $x_1$ and $x_n$ of $F$.
\label{fig:generators}}
\end{figure}

The reduced 
tree pair diagrams associated to the generators $x_0, \ x_1$ and
$x_n$ are pictured in Figure \ref{fig:generators}.  As explained in Lemmas
2.6 and 2.7 of \cite{ct}, the generators $x_0$ and $x_1$ can be viewed in terms
of rotations of rooted binary trees as well.  The generator $x_0$
can be interpreted as a left rotation at the root of the left
tree in the pair, yielding the right tree in the pair. Similarly,
the generator $x_1$ performs a left rotation at the right child
of the root, transforming the left tree in the pair to the
right one. The inverses $x_0^{-1}$ and $x_1^{-1}$ perform  right
rotations at the root node and right child of the root node, respectively.
 
Analogously, right multiplication of an element $w$ given by a possibly
unreduced representative  $(T_1,T_2)$ by the generator $x_0^{-1}$
yields the tree pair diagram $(T_1',T_2')$ in which $T_1'$ differs from $T_1$ by a left rotation at the root node, and $T_2' = T_2$.  We can similarly interpret right multiplication by $x_0, \ x_1^{\pm 1}$ and $x_n^{\pm 1}$.

One complication that may arise when using the geometry of the tree pair diagrams to understand rotation distance is the possibility of requiring unreduced representatives in order to perform the group multiplication.  Since elements of Thompson's
group are equivalence classes of tree pair diagrams, we can
always multiply any group element $w$ by any group generator $g$.  It is possible that we may have to add carets to the reduced tree pair diagram for $w$ in order to carry out this multiplication.  From the standpoint of group theory, the reduced and unreduced tree pair diagrams are interchangeable.  When considering rotation distance, however, we are not allowed to change the number of carets in the starting tree.  Thus certain rotations, corresponding to multiplication by specific generators, may not be permitted when calculating rotation distance.

\begin{figure}
\includegraphics[width=5.5in]{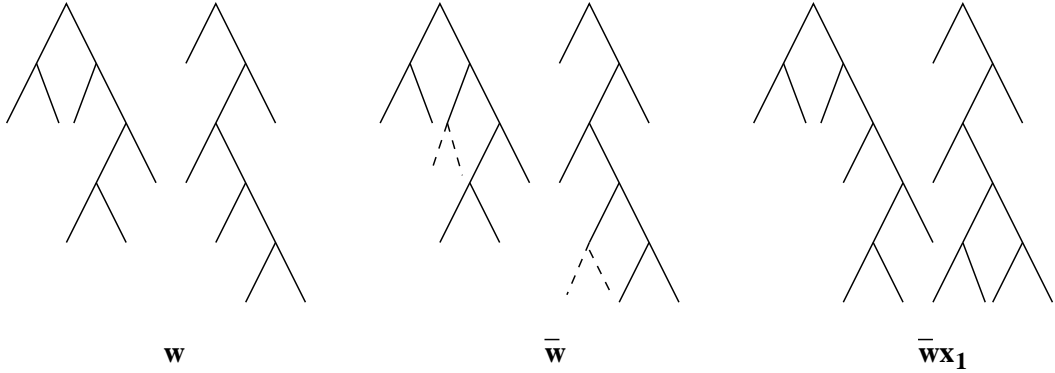} \\
\caption{In order to multiply the tree pair diagram representing $w$ by the generator $x_1$, we form the unreduced representative $\bar{w}$ of $w$ in which the dashed caret is added to both trees.  Then are we able to form the product $\bar{w}x_1$ which is the same group element as $wx_1$.}
\label{fig:addcaret}
\end{figure}

For example, we cannot perform a right
rotation at the right child of the root  to either of the trees in
the tree pair diagram for $x_0$ as shown in Figure \ref{fig:generators}
because neither tree contains a left subtree of the right child of the root.  As an element of Thompson's group, we can enlarge any pair
of trees to another representative in the equivalence class, and thus
are able to multiply any element by any generator.  A typical such application is
shown in Figure \ref{fig:addcaret} where a caret is added to a tree
to be able to perform the desired rotation.  The tree pair diagram $w$ does not have
a left child of the right child of the root, so performing a right rotation at the right child
of the root is not possible.  However, the word $\bar{w}$ which represents the same element
of $F$ does have a left child of the right child of the root and it is possible to perform
the right rotation at the right child of the root there.  We obtain $\bar{w}$ by adding an additional
caret (indicated by dashing) to leaf number 2 in both trees of the tree pair diagram.

To describe when it is necessary to add a caret to a tree to perform
a particular rotation, we make the following definitions.
We say {\em a right rotation at the root can be applied} to a tree $T$
if the left subtree of the root of $T$ is non-empty.  
Similarly, we say  {\em a left rotation at the root can be applied} to a tree $T$
if the right subtree of the root of $T$ is non-empty and we also adopt this terminology
when performing rotations at other nodes along the spine of the tree.

Understanding when rotations can be performed on trees helps
us develop the connection between rotations of trees and right
multiplication by generators of $F$.  If, for example, we have
a tree pair diagram $(T_1, T_1)$ representing the identity and
we can perform a left  rotation at the root to $T_1$ to obtain $ x_0 T_1$,
then the new tree pair diagram $(T_1 , x_0 T_1)$ is the tree
pair diagram representing the word $x_0$ in $F$, and similarly
the  new tree pair diagram $(x_0 T_1, T_1)$ is tree
pair diagram representing the word $x_0^{-1}$ in $F$.

We note that in $F$, multiplication by a generator may result in an unreduced tree pair diagram.  So during the course of
a sequence of multiplications by generators of $F$, the number
of carets in the reduced tree pair diagram representing the partial
products may fluctuate -- rising when it is necessary to add one or more carets to apply a  generator, and
falling when multiplication by a generator results in an unreduced tree pair diagram.
To understand
rotation distance, however, as  we
apply a sequence of rotations to a single tree, we do not allow the
number of carets in the tree to change.

The link between restricted rotation distance
and Thompson's group $F$ is the word metric on $F$ with respect to the generators
$\{x_0,x_1\}$.   Given two rooted binary trees $T_1$ and
$T_2$ with the same number of nodes, we consider a minimal length
word in $x_0^{\pm 1}$ and $x_1^{\pm 1}$ representing the element
$w = (T_1,T_2) \in F$. As described in \cite{rotipl}, this word gives a minimal sequence of
rotations at the root and right child of the root which transform the tree
$T_1$ into the tree
$T_2$. It follows from Fordham \cite{blake:gd} that these minimal words which transform one tree into the other maintain a constant number of carets at each stage in the sequence of rotations.  The issue of certain rotations altering the number of nodes in the tree does not arise in the case of restricted rotation distance.

More precisely, suppose that $w \in F$ is given by the tree pair diagram $(T_1, T_2)$, and a minimal length representative for $w$ is $g_1 g_2 \cdots g_n$, where each $g_i \in \{x_0^{\pm 1},x_1^{\pm 1}\}$.  Then the tree pair diagram $(T_1, g_n \cdots g_2 g_1 T_1)$ will represent $w$, and we can think of the sequence of
generators $g_n \cdots g_2 g_1$ as a sequence of rotations which transforms $T_1$ to $T_2$.
At each stage of this process, we will be able to perform the rotation corresponding to the generator $g_{i+1}$ to the tree $g_i \cdots g_2 g_1 T_1$
without adding extra carets.  There may be reductions possible to tree pair diagrams, or equivalently to the normal forms, during this process, but
from the standpoint of rotation distance we do not want to take advantage of
these reductions.  Instead, we keep the number of carets constant at each
stage.

Equivalently, we can think of $(g_n \cdots g_2 g_1 T_1, T_2)$ as a representative of
the identity and witness the transformation of $T_1$ to $T_2$ by considering the
sequence of tree pair diagrams $$(T_1, T_2), (g_1 T_1, T_2), \ldots,  (g_n \cdots g_2 g_1 T_1, T_2).$$

Below, we consider other possible locations for rotations to occur, and again 
 exploit the link to Thompson's group $F$, but now considering other appropriate
 generating sets for $F$, where the generators are chosen to reflect the locations where rotations
 are permitted.  We assign a level to each node or caret in the tree as follows.  The root node is defined to have level zero.  The level of a node $N$ is the number of edges in a minimal length path connecting $N$ to the root node.  The level of a caret $C$ is defined to be the level of the node associated to that caret.

Writing the generators $x_n$ for $n >1$ via the relators $x_n = x_0^{-(n-1)} x_1 x_0^{n-1}$, we relate each generator to the following rotation of a tree $T$. We denote the all-right tree with the appropriate number of carets by $*$.  Group multiplication must be between a pair of elements, and each element 
corresponds to a pair of trees, so we use the tree $*$ as the positive tree corresponding to $T$.  The product of the generator $x_n$ and the tree pair diagram $(T,*)$ performs a right rotation to $T$ at the caret at level $n$ along the right arm of $T$.  In all that follows, when we describe a generator as inducing a rotation on a single tree $T$ rather than on a tree pair diagram, we are forming the product with the pair $(T,*)$ as above.

\section{Metrics on $F$ and rotation distances}

\subsection{Relation to the word metric.}

In \cite{rotbound}, the word length with respect to the finite
generating set $\{x_0, x_1\}$ of $F$ is used to compute the
restricted rotation distance between a pair of trees, using techniques of Fordham
\cite{blake:gd}.  Fordham developed a method for computing the exact length of
an element of $F$ directly from the reduced tree pair diagram representing that element.


\begin{definition}
If $T_1$ and $T_2$ are trees with the same number of nodes,
we define the {\em restricted rotation distance} $d_{RR}(T_1,T_2)$
as the minimal number of rotations required to transform $T_1$ to $T_2$,
where rotations are allowed at the root and the right child of the root.
\end{definition}

Restricted rotation distance is well-defined for any two trees with
the same number of leaves, as shown in \cite{rotipl}. We then obtain
the following sharp bound on restricted rotation distance.

\begin{theorem}[\cite{rotbound}, Theorems 2 and 3]
\label{thm:rotbound} Given two rooted binary trees $T_1$ and $T_2$
each with $n$ nodes, for $n \geq 3$, the restricted rotation
distance 
between them satisfies $d_{RR}(T_1,T_2) \leq 4 n- 8$.
Furthermore, for $n\geq 3$, there are trees $T_1'$ and $T_2'$ with $n$ nodes realizing this
bound;  that is, with  $d_{RR}(T_1',T_2') = 4 n- 8$.
\end{theorem}

Intermediate between the two-element generating set $\{x_0,x_1\}$ and the infinite
generating set $\{x_0, x_1, \ldots \}$ are other finite generating sets of the form $\{x_0,x_{i_1}, \ldots ,x_{i_L}\}$, where we arrange the indices of the (distinct) generators in increasing order.  Analyzing the infinite generating set corresponds to allowing all rotations along the right side of the tree.  Finite generating sets correspond to allowing finite collections of rotations at the root node and other nodes along the right side of the tree.

\begin{definition}
Let ${\mathcal S}=\{x_0,x_{i_1}, \ldots ,x_{i_L}\}$ be a finite subset of the infinite generating set for $F$
and $T_1$ and $T_2$ be trees with the same number of leaves.
We define  $d_{RRA}^{\mathcal S} (T_1,T_2)$,  the {\em  restricted right-arm rotation distance
with respect to  ${\mathcal S}$}, as the minimal number of rotations required to transform $T_1$ to
$T_2$, where the rotations are only allowed at levels $0, i_1, \ldots, i_{L-1}$ and
$i_L$  along
the right side of the tree.
\end{definition}

We will see below that unlike restricted rotation distance, restricted right-arm rotation distance may not be defined between all pairs of trees with the same number of nodes.
We use the notation $| \cdot |_{\mathcal S}$ to denote the word length of an element of $F$ with respect to the generating set $\mathcal S$.  
We now relate the restricted right-arm rotation distance $d^{\mathcal S}_{RRA}(T_1,T_2)$ to $| (T_1,T_2) |_{\mathcal S}$.

We consider two trees $T_1$ and $T_2$ each with $n$ nodes.  The word length of the element $w=(T_1,T_2) \in F$ with respect to a generating
set ${\mathcal S}$ is the length of the shortest expression for $w$ in that generating set.   However, when considering the corresponding rotations to the tree pair diagram for $w$,
 we have no analogue of Fordham's proof that a minimal length representative
 in these generators can be constructed while maintaining a constant number of nodes in each tree. Thus, it may be possible that a minimal length representative for $w=(T_1,T_2) \in F$ with respect to ${\mathcal S}$ includes some rotations which would require the addition of carets to the trees and are thus
 not permitted. 
Therefore, we see that the word length $|(T_1,T_2)|_{\mathcal S}$ provides only a lower bound on the rotation distance $d_{RRA}^{\mathcal S}(T_1,T_2)$, when this rotation distance is defined. If this word length corresponds to a sequence of rotations in which the number of nodes remains constant at each intermediate step, then we have computed the actual restricted right-arm rotation distance between the two trees.  These cases will be addressed below.

\begin{figure}
\includegraphics[width=3in]{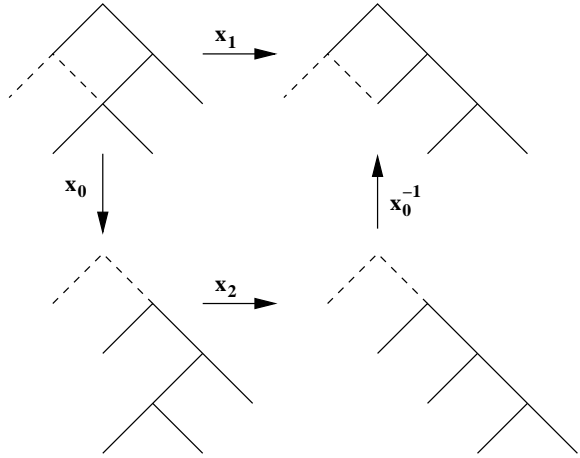} \\
\caption{A right rotation at the right child of the root performed on the reduced tree by $x_1$ and on the partially reduced tree by $x_0x_2x_0^{-1}$.  The trees $T_1$ and $T_2$ are the solid trees
in the left and right respectively of the top row, and the trees  $T_1'$ and $T_2'$ are respectively the
same trees with the additional single dashed caret in each case.
\label{fig:partialred}}
\end{figure}

For example, we consider the trees shown in Figure \ref{fig:partialred}.  The desired transformation from the top left tree $T_1$ drawn in solid lines  to  the top right tree $T_2$ drawn in solid lines would be given by $x_1$, a single right rotation at the right child of the root.  But if the permitted locations for rotation are only at the root (corresponding to the generator $x_0^{\pm1}$)  and the right child of the right child of the root (corresponding to the generator $x_2^{\pm1}$),  it will be impossible to accomplish the desired transformation without adding extra nodes, and the corresponding restricted right-arm rotation distance is not defined between those two trees.

If we are permitted to add a node to the leftmost leaf of each tree, as shown with the dashed carets, to obtain the related problem of transforming the new tree $T_1'$ into $T_2'$ (drawn including the dashed carets) then the transformation  would be possible using only the allowed rotations.  The unreduced form of the top tree pair diagram $(T_1',T_2')$ drawn including the dashed caret  is $x_0 x_2 x_0^{-1}$ which reduces to $x_1$ in the usual manner in $F$, if desired.   If rotations are permitted at the root and right child of the root, the rotations that transform  $T_1$ to $T_2$ are exactly the same as those to perform the transformation from $T_1'$ to $T_2'$ and the added dashed caret is simply carried along intact.  However, if we are only permitted to rotate at the root and right child of the right child of the root, the added caret is essential in allowing that transformation, though it does take two additional steps.  We cannot transform  $T_1$ to $T_2$ but we can easily transform $T_1'$ to $T_2'$ by rotating rightwards at the root, rightwards at the right child of the right child of the root, and then leftwards at the root, as pictured.

We can describe exactly when a tree $T_1$ can be transformed into $T_2$ without adding nodes with respect to a specified set of allowed rotations along the right-arm of the tree; that is, exactly
when the restricted right-arm rotation distance is defined. First, we consider the case when the word in normal form associated to $(T_1,T_2)$ is already reduced; that is, when $(T_1, T_2)$ is a reduced tree pair diagram.
 
\begin{lemma}\label{lemma:reducednotdefined}
Let ${\mathcal S}=\{x_0,x_{i_1}, \ldots ,x_{i_L}\}$ be a generating set for $F$
 with $0<i_1 < i_2 < \cdots < i_L$.   We
 consider the corresponding restricted right-arm rotation distance $d_{RRA}^{\mathcal S}$,
 where rotations are allowed at nodes at levels $0, i_1, \ldots, i_L$ on the right side of the tree.
Suppose $T_1$ and $T_2$ are finite rooted binary trees with the same number of nodes forming a
reduced tree pair diagram $w = (T_1,T_2) \in F$ with unique normal form given by $$x_{i_1}^{r_1} x_{i_2}^{r_2}\cdots x_{i_k}^{r_k}
x_{j_l}^{-s_l} \cdots x_{j_2}^{-s_2} x_{j_1}^{-s_1}.$$
If $x_t^{\pm 1}$ for $1 \leq t \leq i_1-1$ appears in this
normal form, then the restricted right-arm
rotation distance $d_{RRA}^{\mathcal S}(T_1,T_2)$ is not defined.  
Conversely, if no $x_t^{\pm 1}$ for $1 \leq t \leq i_1-1$ appears in the unique
normal form, then the restricted right-arm
rotation distance $d_{RRA}^{\mathcal S}(T_1,T_2)$ is defined.  
\end{lemma}

When $i_1=1$, it follows from Lemma \ref{lemma:reducednotdefined} that $d_{RRA}^{\mathcal S}$ will always be defined.  This includes the special case of restricted rotation distance, when ${\mathcal S} = \{x_0,x_1\}$.

\begin{proof}
We recall that the leaf exponent of the leaf numbered $n$ in a tree is the length of the maximal path of left edges from leaf $n$ which does not reach the right side of the tree.  Observe that the leaf exponent that changes as a result of a rotation at the node at level $h$ on the right arm of the tree corresponds to the leftmost leaf in the left subtree of the node where the rotation occurs.

First, we suppose that  $x_t^{\pm 1}$ for $1 \leq t \leq i_1-1$ appears in the unique
normal form for $(T_1,T_2)$. So $t$ appears as the leaf number of a left leaf of a caret in either in $T_1$ or $T_2$ or possibly both.
If the restricted right-arm rotation distance  $d_{RRA}^{\mathcal S}(T_1,T_2)$ is defined,
then the sequence of rotations transforming $T_1$ into $T_2$ does not change the number of nodes in the tree at any intermediate step and thus no leaves are added or removed during this process.
So there is no potential renumbering of leaves, as there may be when considering
equivalence classes of tree pairs in Thompson's group $F$.
We consider the leaf numbers whose exponents can be affected by rotations at the permitted nodes. Rotations are permitted at the root node and at levels $i_j$ along the right side of the tree.  Rotations at the root can affect only the exponent of leaf zero, as it will be the leftmost leaf in the left subtree attached at the root.   Other rotations can affect the exponents of leaves which are the leftmost leaves of left subtrees of right nodes at levels $i_1$ and lower.  The left subtree of the right node at level $h$ will have leaves numbered at least $h$, so if $t<i_1$, then no rotation at level $i_j$ can affect the exponent of leaf $t$.  So if there is a left leaf in the range $1 \leq t \leq i_1-1$ present in $T_1$, rotations at the root cannot affect its exponent, and rotations at levels $i_1$ and greater cannot affect its exponent. 


If leaf $t$ has different exponents in $T_1$ and  $T_2$, since the allowed rotations cannot change its exponent, $T_1$ cannot be
transformed into tree $T_2$ by the allowed rotations.  If leaf $t$ is present in both
trees with the same exponent, then since $w$ is in unique normal form, the
exponent of leaf $t+1$ must also be non-zero in at least one of the trees.  Moreover, leaves numbered $t$ and $t+1$ belong to the left subtree of the same node on the spine.  Thus none of the allowed rotations can affect the leaf exponent of leaf $t+1$ as well.  We iterate this argument with leaves $t+1$ and $t+2$.  Thus, we see that if any   $x_t^{\pm 1}$  with $1 \leq t \leq i_1-1$ appears, then the two trees cannot be connected by any sequence of the allowed rotations without the addition of extra nodes.

Conversely,  if  $x_t^{\pm 1}$ for $1 \leq t \leq i_1-1$ do not appear in the normal form, then we can rotate $T_1$ rightwards at the root  by application  of an appropriate power of $x_0^k$ so that all of the nontrivial subtrees then hang from the right arm of the tree at levels $i_1$ and greater. In Proposition \ref{prop:infgenset} we show that the right-arm rotation distance is always defined between two trees with the same number of nodes.  This allows us to finish the proof with the following argument.

We now use $x_0, x_{i_j}$ and conjugates of   $x_{i_j}$ by powers of $x_0$  to rotate the tree to an all-right tree, just as in the infinite generating set, without adding any extra nodes.  So we can
transform $T_1$ to the all-right tree, and then from the all-right tree,
we can again use $x_0, x_{i_1}$ and conjugates of   $x_{i_1}$ by powers of $x_0$ 
(and possibly other $x_{i_j}$, if desired) to 
transform the all-right tree to $T_2$ without adding extra nodes.  Thus,  $d_{RRA}^{\mathcal S}(T_1,T_2)$ is defined.  There may be more efficient ways of accomplishing this transformation but it is clear that there is at least one way of doing it without adding extra nodes,
so the restricted right-arm rotation distance is defined.
\end{proof}

To understand the case where $(T_1,T_2)$ is an unreduced tree pair diagram, and thus we do not obtain the unique normal for the element directly from the leaf exponents,
  we introduce the notion of {\em partial reduction}.
Partial reduction is similar to ordinary reduction except that we do not want to remove left nodes common to both trees.  Stated algebraically, it means that if the normal form for the element contains instances of $x_0$ and $x_0^{-1}$ but not $x_1^{\pm 1}$, we do not simplify the expression, as we do when $x_k$ and $x_k^{-1}$ appear but not $x_{k+1}^{\pm 1}$ for $k > 0$.  The presence of these additional left nodes may allow us to perform rotations which would not be permitted otherwise without increasing the number of carets in the trees.  This phenomenon occurs in the tree pairs
shown in Figure \ref{fig:partialred}, where the restricted right-arm rotation distance where rotations are
permitted at the root and the right child of the right child of the root is defined between
$T_1'$ and $T_2'$ but not between $T_1$ and $T_2$.

\begin{definition}
A word $w$ in $F$ in normal form  is {\em partially reduced} if it is of the form
$x_{i_1}^{r_1} x_{i_2}^{r_2}\cdots x_{i_k}^{r_k} x_{j_l}^{-s_l}
\cdots x_{j_2}^{-s_2} x_{j_1}^{-s_1} $ with $0 \leq i_1<\cdots < i_k$ and
$0 \leq j_1< \cdots < j_l$, with $r_n$ and $s_n$ all positive, and 
if we further require the partial
reduction condition that for $i>0$, when both $x_i$ and $x_i^{-1}$ occur, so
does at least one of $x_{i+1}$ or  $x_{i+1}^{-1}$.
\end{definition}

For any word $w$ in (not necessarily unique) normal form, there will be a maximal length word $w'$ satisfying
the partial reduction condition which we can easily obtain using the procedure described above.

The partial reduction allows us to prove the following lemma, which describes when one given tree can be transformed into
another with respect to a specified set of rotations, when the initial tree pair diagram is unreduced.   The proof is identical to that of Lemma \ref{lemma:reducednotdefined}.

\begin{lemma}\label{lemma:notreducednotdefined}
Let ${\mathcal S}=\{x_0,x_{i_1}, \ldots ,x_{i_L}\}$ be a generating set for $F$
 with $0<i_1 < i_2 < \cdots < i_L$.   We
 consider the corresponding restricted right-arm rotation distance $d_{RRA}^{\mathcal S}$
 where rotations are allowed at the root node and at right nodes of levels $i_1, \ldots, i_L$.
  Suppose $T_1$ and $T_2$ are finite rooted binary trees with the same number of nodes forming a
 tree pair diagram $w = (T_1,T_2) \in F$ and  that $w$ has the partially reduced normal form of
 maximum length given by the word $$w'= x_{i_1}^{r_1} x_{i_2}^{r_2}\cdots x_{i_k}^{r_k}
x_{j_l}^{-s_l} \cdots x_{j_2}^{-s_2} x_{j_1}^{-s_1}.$$
If $x_t^{\pm 1}$ for $1 \leq t \leq i_1-1$ appears in this partially
reduced normal form, then the restricted right-arm
rotation distance $d_{RRA}^{\mathcal S}(T_1,T_2)$ is not defined.  
Conversely, if no $x_t^{\pm 1}$ for $1 \leq t \leq i_1-1$ appears in this partially
reduced normal form, then the restricted right-arm
rotation distance $d_{RRA}^{\mathcal S}(T_1,T_2)$ is defined.  
\end{lemma}

When rotations at all nodes along the right side of the tree are allowed, we obtain the right-arm rotation distance $d_{RA}$, understood by the methods of   Culik and Wood \cite{cw} and Pallo \cite{palloright}.  Culik and Wood considered general rotation distance but used only rotations on the right
side of the tree to show their upper bound of $2n-2$, while Pallo intentionally restricts to only allow rotations on the right hand side of the tree.
Pallo's situation is analogous to restricted rotation distance, which considers only the rotations corresponding to the generators $x_0$ and $x_1$, because the word length once again yields the exact rotation distance.

\begin{proposition}\label{prop:infgenset}
Let ${\mathcal I}$ denote the standard infinite generating set for $F$, and let $T_1$ and $T_2$ be binary trees, each with $n$ nodes. Then
$$d_{RA}(T_1,T_2) = |(T_1,T_2)|_{\mathcal I}.$$
\end{proposition}

\begin{proof}
We will assume that the tree pair diagram $(T_1,T_2)$ is reduced.  If it is not, we form the tree pair diagram $(T_1',T_2')$ representing the same group element which is reduced.  The rotations necessary to transform $T_1'$ into $T_2'$ will also transform $T_1$ into $T_2$, since no additional rotations are necessary to alter the nodes which cause $T_1$ and $T_2$ to be unreduced.  The
nodes which were removed during the reduction are 
identical in both trees and are carried along unchanged during the rotations which transform $T_1'$ to $T_2'$.
The leaf exponent method of associating the unique normal form to the tree pair diagram 
described above shows that each tree provides one part of the normal form; in the pair $(T_1,T_2)$ the tree $T_1$ corresponds to the terms with negative exponents and $T_2$ to those with positive exponents.  We thus write the normal form as the product $P N$, where $N$ contains the generators with negative exponents, and $P$ those with positive exponents.

We see that $N$ is a word which rotates the tree $T_1$ into the all-right tree 
without requiring the addition of any nodes, and  the subword $P$ is a string of
generators which rotates the all-right tree into the tree $T_2$.  

Thus we see that a lower bound for right arm rotation distance is $|(T_1,T_2)|_{\mathcal I}$, and an upper bound is given by combining the length of the strings $P$ and $N$.  It follows from Lemma \ref{lemma:geodesic} that $|(T_1,T_2)|_{\mathcal I} = |(T_1,*)|_{\mathcal I} + |(T_2,*)|_{\mathcal I}$ where $*$ is the all-right tree with $n$ nodes, proving the proposition.
\end{proof}

\subsection{Bounds on restricted rotation distances}

Now that we have described the relationship between the different rotation distances and word lengths in $F$, we obtain numerical bounds on these rotation distances as summarized in Table \ref{summarytable}.
We note that word length of an element of $F$ computed with respect to a generating set of the form ${\mathcal S}$ given above has the potential to be much shorter than the word length of the same element computed with respect to the generating set $\{x_0,x_1\}$.  Thus we might expect significantly smaller asymptotic upper bounds on restricted right-arm rotation distance than on restricted rotation distance.  In fact, this is not the case, and the difference between the upper bounds on the two rotation distances is at most a constant.

The goal of this section is to prove that the multiplicative constant of $4$ in the upper bound on restricted rotation distance cannot be improved upon when we allow additional rotations along the right arm of the tree.  Both of these rotation distances between two trees with $n$ nodes each, when defined, are bounded above by $4n$ minus a constant.  This constant depends upon the particular finite set of rotations permitted.  These bounds are shown to be sharp for restricted rotation distance in \cite{rotbound}.  We show below that they are asymptotically sharp for restricted right-arm rotation distance as well.  While allowing additional rotations may shorten the restricted right-arm rotation distance between certain pairs of trees, asymptotically the worst-case scenario differs from restricted rotation distance only by an additive constant.  One way to improve the multiplicative constant of $4$ is to allow rotation at an infinite collection of nodes along the
right side of the tree, in which case the multiplicative constant of the bound may decrease to $2$.

The necessity of the constant $4$ is shown in two steps.  We first show that the restricted right-arm rotation distance, when defined, is always bounded above by $4n-8$, where $n$ is the number of nodes in either tree.  We then show that there are words which approach this bound to within an additive constant.

\begin{proposition}\label{prop:rrabounds}
Let ${\mathcal S} = \{x_0, x_{i_1}, x_{i_2}, \ldots, x_{i_L}\}$ be a generating set for $F$ with $0<i_1 < i_2 < \cdots < i_L$, and let $d_{RRA}^{\mathcal S}$ be the corresponding restricted right-arm rotation distance.  Let $T_1$ and $T_2$ be binary trees, each with $n$ nodes with $n \geq 3$, for which $d_{RRA}^{\mathcal S}(T_1,T_2)$ is defined.  Then 
$$
d_{RRA}^{\mathcal S}(T_1,T_2) \leq 4n-8.
$$
\end{proposition}

\begin{proof}
The case where $i_1=1$ is already addressed by the analysis of ordinary rotation distance,
described in \cite{rotbound}.
We consider the element $w=(T_1,T_2) \in F$, where $T_1$ and $T_2$ are trees for which
the relevant restricted right-arm rotation distance $d_{RRA}^{\mathcal S}$ is defined, and assume that $i_1 > 1$.

{\bf Case 1: The tree pair diagram $(T_1,T_2)$ is reduced.} 

In this case,  we know that the normal form of $w$ contains no generators $x_t^{\pm 1}$ for $1 \leq t \leq i_1-1$.  In addition, this normal form can contain $x_0$ or $x_0^{-1}$ but not both.  If both $x_0$ and $x_0^{-1}$ were present in the normal form with no $x_1^{\pm 1}$ generator, then the normal form could be reduced.  We can  assume by symmetry 
that the normal form for $w$ contains $x_0^{-k}$ but no factors of $x_0$. 

Using the correspondence between the normal form and the leaf exponents in the trees $T_1$ and $T_2$, we see that the leaves of both trees numbered from 1 through $i_1-1$ are either exposed right leaves of left nodes or exposed left leaves of right nodes.  In $T_1$, denote the (possibly empty) subtrees of the left and right nodes by $A_1,A_2, \ldots ,A_n$, where the smallest leaf number in $A_1$ is $i_1$.  If $A_1$ is empty, by "smallest leaf number", we mean the number of the leaf attached
to the spine of the tree in that position.
  Similarly, in $T_2$ denote these subtrees by $B_1,B_2, \ldots ,B_m$, where the smallest leaf number in $B_1$ is $i_1$.  

\begin{figure}
\includegraphics[width=4in]{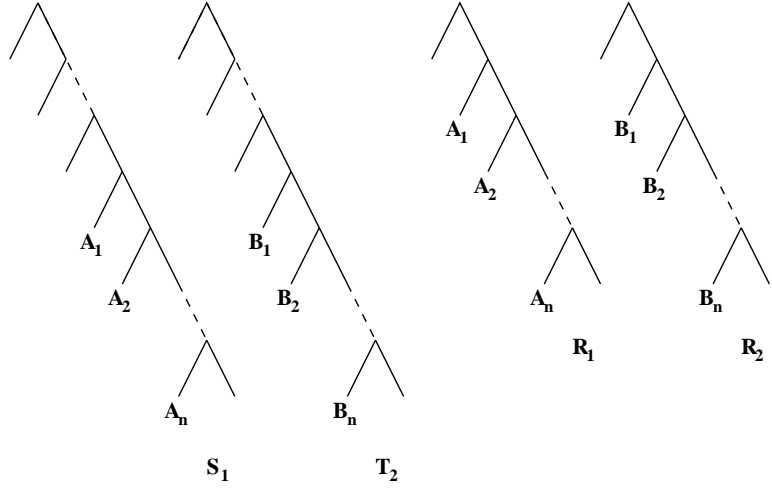}\\
\caption{The tree pair diagrams for the words used to show that the restricted right-arm rotation distance is bounded above by $4n$.}
\label{fig:T1T2}
\end{figure}

Let $w' = wx_0^{k}$, so that the tree pair diagram $(S_1,T_2)$ of $w'$ has tree $S_1$ containing a single left node, namely the root node, and $i_1-1$ right nodes with exposed left leaves, followed by right nodes having $A_1, \ldots ,A_n$ as their left subtrees.  The pair $(S_1,T_2)$ has the form given in Figure \ref{fig:T1T2}.

We consider the element $v \in F$ which has tree pair diagram $(R_1,R_2)$, where $R_1$ has a single left node, namely the root node, and the left subtree of the right node at height $i$ is $A_i$.  The tree $R_2$ is defined analogously, using the subtrees $B_i$ from the original tree $T_2$.  Since restricted rotation distance is well defined for all trees with the same number of nodes, we apply Theorem \ref{thm:rotbound} to obtain the bound $d_{RR}(R_1,R_2) \leq 4(n-(i_1-1)) - 8$.  This restricted rotation distance is realized by a string $\alpha$ of the generators $\{x_0^{\pm 1},x_1^{\pm 1}\}$.

We define a string of generators $\alpha'$ by replacing each instance of $x_1^{\pm1}$ in $\alpha$ with $x_{i_1}^{\pm1}$.  Then this string of generators exactly produces the tree pair diagram $(S_1,T_2)$.  Since the number of nodes in each tree remains constant as each generator from $\alpha$ is applied to create $(R_1,R_2)$, the same is true as we multiply the generators in $\alpha'$ to create $w'=(S_1,T_2)$.

Thus the restricted right-arm rotation distance with respect to ${\mathcal T} = \{x_0,x_{i_1}\}$ is bounded as follows:
$$d_{RRA}^{\mathcal T}(S_1,T_2) \leq 4(n-(i_1-1)) - 8.$$ 

Now we note that $w = w'x_0^{-k}$, and since there were initially $k+1$ left nodes in the tree $T_1$, the number of nodes in each tree remains constant during these successive multiplications by $x_0^{-1}$.  Thus the string $\alpha'x_0^{-k}$ realizes the restricted rotation distance between the trees $T_1$ and $T_2$.

If $k \leq i_1-1$, then the left nodes which are changed to right nodes under multiplication by $x_0^{-k}$ do not appear in $R_1$ and $R_2$, and so are not represented in the upper bound given above.  Thus, when the rotation distance is increased by $k$, we trivially extend the bound to 
$$d_{RRA}^{\mathcal T}(T_1,T_2) \leq 4n - 8.$$
Since adding extra generators to the generating set, or equivalently allowing rotations at additional nodes, can only decrease the rotation distance, the upper bound still holds when we consider the entire generating set ${\mathcal S}$.

If $k \geq i_1$, then the left nodes which are changed to right nodes by these multiplications by $x_0^{-1}$ are of two types: those with exposed left leaves numbered from $1$ to $i_1-1$, and those with left subtrees of the form $A_i$.  The first type of right node is not counted in the upper bound given above, and thus we increase the number of nodes in the bound by $i_1-1$ to (more than) account for the additional generators.  

The right nodes of the second type, with left subtrees of the form $A_i$, are already counted in the bound given above.  However, we recall that the word $\alpha'$ which realizes the restricted right-arm rotation distance between $S_1$ and $T_2$, came from the word $\alpha$ in $\{x_0^{\pm 1},x_1^{\pm 1}\}$.  We know from Fordham's method of calculating word length with respect to the generating set $\{x_0,x_1\}$ directly from the tree pair diagram that each pair of nodes with the same infix number in each tree contributes a certain number of generators to this word length.  Fordham calls this the {\em weight} of the pair of nodes.  We see from Fordham's table of weights \cite{blake:gd} that any pair of nodes in which one node is a right node has a weight of at most three.  So using an extra generator of the form $x_0^{-1}$ to transform this right node into a left node means that these nodes contribute at most four generators each to the length of the word realizing the restricted right-arm rotation distance between $T_1$ and $T_2$.  We have thus shown the existence of the upper bound
$$d_{RRA}^{\mathcal T}(T_1,T_2) \leq 4n - 8.$$
Since rotation distance can only decrease when additional rotations are permitted,  this extends immediately to show 
$$d_{RRA}^{\mathcal S}(T_1,T_2) \leq 4n - 8.$$

{\bf Case 2: The tree pair diagram $(T_1,T_2)$ is not reduced.}
 
In this case, since $d_{RRA}^{\mathcal S}(T_1,T_2)$ is defined,
we know from Lemma \ref{lemma:notreducednotdefined} that the unreduced tree pair diagram $(T_1,T_2)$ yields a word $\alpha$ which is partially reduced and represents $w$.  We obtain $\alpha$ by considering the unreduced normal form arising from $(T_1,T_2)$ and applying the usual reduction rules but without reducing instances of $x_0$ and $x_0^{-1}$ with no $x_1^{\pm 1}$.  We may also be able to partially reduce $(T_1,T_2)$ to correspond to this partially reduced normal form for $w$; in this case the number of nodes in the tree pair diagram may reduce to $n' < n$.  The proof of this case is now identical to that of Case 1.  This produces an upper bound of $4n'-8 < 4n-8$ on the restricted right-arm rotation distance between the two trees.
\end{proof}

We now show that the multiplicative constant of 4 is necessary for the above inequality.

\begin{theorem}\label{thm:rra-sharp}
Let ${\mathcal S}=\{x_0,x_{i_1}, \ldots ,x_{i_L}\}$.  Then there exist trees $T_1$ and $T_2$, each  
 with $n$  nodes, so that $d_{RRA}^{\mathcal S}(T_1,T_2)$ is defined, and with 
 $$d_{RRA}^{\mathcal S}(T_1,T_2) \geq 4n-4 i_L - 4.$$
\end{theorem}

The generating set ${\mathcal S}$ used in Theorem \ref{thm:rra-sharp} corresponds to a series of rotations along the right side of the tree from levels $0$ to $i_L$ but does not necessarily
 include all rotations at levels within this range.  We now enlarge our generating set to correspond to all rotations at levels $0$ to $i_L$, and work with this set ${\mathcal S}'$ in Theorem \ref{thm:rra-sharpfull}.  It will be enough to use this larger set of generators and show that  $d_{RRA}^{\mathcal S'}(T_1,T_2) \geq 4n-4 i_L -4$.
Thus we prove the following theorem.

\begin{theorem}\label{thm:rra-sharpfull}
Let ${\mathcal S'}=\{x_0,x_{1},x_2, \ldots ,x_{i_L}\}$, with $m=i_L$.  Then there exist trees $T_1$ and $T_2$, each with $n$  nodes, so that $d_{RRA}^{\mathcal S'}(T_1,T_2)$ is defined, and with $$
d_{RRA}^{\mathcal S'}(T_1,T_2) \geq 4n-4 m - 4.$$
\end{theorem}

\begin{figure}\includegraphics[width=4in]{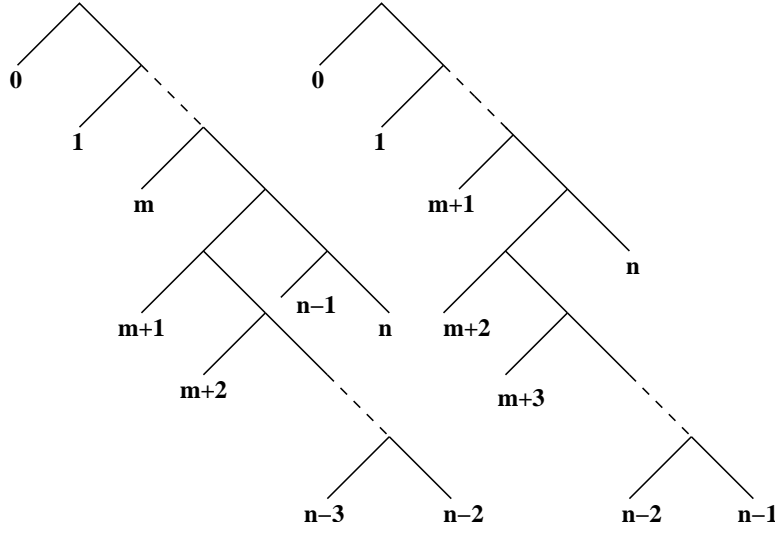}\\ 
\caption{The tree pair diagram $(T_1,T_2)$ for words of the form $x_{m+2} x_{m+3} \cdots x_{n-2} x_{n-3}^{-1} x_{n-4}^{-1}\cdots x_{m+1}^{-1}$.}
\label{badword}
\end{figure}

The elements we will use to prove this theorem have normal form
$$x_{m+2} x_{m+3} \cdots x_{n-2} x_{n-3}^{-1} x_{n-4}^{-1}\cdots x_{m+1}^{-1}$$ and tree pair diagram which we denote $(T_1,T_2)$.  These elements are shown in Figure \ref{badword}.
 
Fordham's method for computing exact word length is only valid for the
generating set $\{x_0,x_1\}$, so we bound the lengths of these elements indirectly in a series of lemmas
which analyze how many generators are needed to change a ``deeply buried'' part of the tree. 

\begin{figure}
\includegraphics[width=3in]{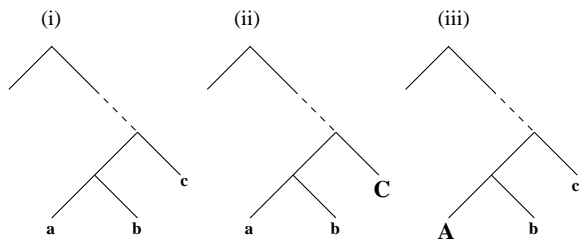}
\caption{Instances where right rotation at node $N$ along the right side of the tree creates or destroys sibling pairs in the tree $T$.  Capital letters represent nonempty subtrees of $T$, and lower case letters denote leaf numbers.
\label{fig:sibpairright}}
\end{figure}

We write $[i,i+1]$ if leaves $i$ and $i+1$ form a sibling pair.  Performing a rotation corresponding to the generator $x_n$ on a tree $T$ is equivalent to taking the product of $x_n$ with the tree pair diagram $(T,*)$, where $*$ is the tree consisting only of the root node and a series of right nodes.  By analyzing the effect of a rotation at a right node on a tree $T$, we see that there are only three configurations of $T$ which allow a sibling pair to be created or destroyed.  These are presented in Figure \ref{fig:sibpairright}, where capital letters refer to nonempty subtrees of $T$ and lower case letters denote leaf numbers.  Right rotation at the appropriate node $N$ along the right side of the tree has the following effect on the sibling pairs.
\begin{itemize}
\item[(i)] The pair $[a,b]$ is destroyed and the pair $[b,c]$ is created.
\item[(ii)] The pair $[a,b]$ is destroyed.
\item[(iii)] The pair $[b,c]$ is created.
\end{itemize}

\begin{figure}
\includegraphics[width=3in]{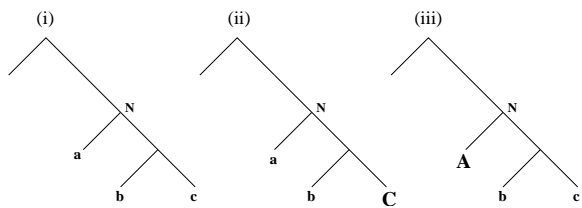}
\caption{Instances where left rotation at node $N$ along the right side of the tree creates or destroys sibling pairs in the tree $T$.  Capital letters represent nonempty subtrees of $T$, and lower case letters denote leaf numbers.
\label{fig:sibpairleft}}
\end{figure}

We can similarly consider left rotation at the node $N$, in which case we refer to Figure \ref{fig:sibpairleft}.  Left rotation at node $N$ along the right side of the tree has the following effect on the sibling pairs.
\begin{itemize}
\item[(i)] The pair $[a,b]$ is created and the pair $[b,c]$ is destroyed.
\item[(ii)] The pair $[a,b]$ is created.
\item[(iii)] The pair $[b,c]$ is destroyed.
\end{itemize}

From these observations, we can see immediately at which nodes it is possible to create
and destroy sibling pairs with a set of rotations.

\begin{lemma}\label{affectedsibs}
Suppose a tree $T'$ is obtained from a tree $T$ by applying a right rotation at a node $N$ at level $n$ on the
right side of $T$. 
\begin{enumerate}
\item If leaves $m$ and $m+1$ are siblings in $T$ and are not siblings in $T'$, then
leaves $m$ and  $m+1$ are the leaves of an exposed node whose parent is node $N$.   
\item If we have the sibling pair $[m,m+1]$ in $T'$ but
not in $T$, then $[m,m+1]$ must be the rightmost node in $T'$ and $m$ must be a leaf in $T$ whose parent is the node at level $n$ in $T$.  

\end{enumerate}
Similarly, the opposite conditions hold for left rotations.
\end{lemma}

We note that when $N$ is the root node, the only sibling pairs that might be affected by rotation at $N$ consist of the first two and the last two leaves in the tree.

When we consider the trees $T_1$ and $T_2$ in Figure \ref{badword}  corresponding to the reduced tree pair diagram representing the element 
$x_{m+2} x_{m+3} \cdots x_{n-2} x_{n-3}^{-1} x_{n-4}^{-1}\cdots x_{m+1}^{-1}$, we see that in $T_1$, leaves
$n-3$ and $n-2$ are siblings and in $T_2$, leaves $n-2$ and $n-1$ are siblings. 

Now we consider the minimal number of rotations needed to change the
sibling pairings from $[n-3,n-2]$ to $[n-2,n-1]$, expressed as a word $w=g_1 g_2 \cdots g_l$, where each $g_i \in {\mathcal S'}$. We will need to first destroy
the sibling pair $[n-3,n-2]$ and subsequently create the sibling pair $[n-2,n-1]$.
The exposed nodes with siblings $[n-3,n-2]$ and $[n-2,n-1]$ are ``deeply buried" in the sense that many rotations are required to affect those nodes and thus those leaf pairings.  We measure this depth more precisely with the following definition.

\begin{definition}
Let $c$ be an exposed caret, and $\alpha_c$ the minimal path from the node of $c$ to the spine of the tree.  The node which is the endpoint of $\alpha_c$ lying on the spine of the tree is called the {\em spinal ancestor} of $c$.  

Define $G(c) = (r,s)$ where $r$ is the number of edges in the path $\alpha_c$ and $s$ is the level in the tree of the spinal ancestor of $c$.
\end{definition}

Note that the spinal ancestor of $c$ can be either a left or right caret.
For example, in the tree $T_1$ for $w$ given in Figure \ref{badword}, we  consider
$c$ as the caret with exposed leaves $[n-3,n-2]$.  The length of the path $\alpha_c$ is $n-m-3$, so $G(c) = (n-m-3,m+1)$.  If $c'$ is the caret in $T_2$ with exposed leaves $[n-2,n-1]$, then $G(c) = (n-m-3,m+2)$. 

In Table \ref{tab:right-G(c)} below, we summarize the changes in $G(c) = (r,s)$ when a single rotation at level $k>0$ is performed on the right arm of the tree containing the exposed caret $c$.  We label each non-spinal node along the path $\alpha_c$, beginning with the one closest to the spine, as follows. We give the node the label $R$ if  it belongs to a caret which is the right child of its parent, and $L$ if that node belongs to a caret which is the left child of its parent.  If the first spinal ancestor of $c$ is on the right side of the tree, then $\alpha_c$ must begin with the label $L$.  Analogously, if the first spinal ancestor of $c$ is on the left side of the tree, then $\alpha_c$ must begin with the label $R$.

The change in $G(c)$ under a single rotation is governed by two factors:
\begin{enumerate}
\item the relative positions of the levels $k$ and $s$, and
\item the first {\em two} labels along the path $\alpha_c$.
\end{enumerate}
The following tables summarize the changes in $G(c)=(r,s)$ when different rotations are performed at level $k$ along the right side of the tree, so we are assuming that the spinal ancestor of $c$ lies on the right arm of the tree.  If the spinal ancestor of $c$ lies on the left arm of the tree, then $G(c)$ is unaffected by a rotation along the right arm of the tree.

\begin{table}[h]\label{G(c)summary-right}
\begin{tabular}{|c|c|c|c|}
\hline
Direction & Initial labels  & Relative position & Change in \\
 of Rotation & of $\alpha_c$ &  of $k$ and $s$ & $G(c) = (r,s)$  \\
\hline 
 & & & \\
Left rotation & LL or LR & $k < s-1$ & $(r,s-1)$ \\
Right rotation & LL or LR & $k < s-1$ & $(r,s+1)$ \\ 
 & & & \\
\hline
 & & & \\
Left rotation & LL or LR & $k = s-1$ & $(r+1,s-1)$ \\
Right rotation & LL or LR & $k = s-1$ & $(r,s+1)$ \\
 & & & \\
\hline
 & & & \\
Left rotation & LL or LR & $k = s$ & $(r+1,s)$ \\
Right rotation & LL & $k = s$ & $(r-1,s)$ \\
Right rotation & LR & $k = s$ & $(r-1,s+1)$ \\
 & & & \\
\hline
 & & & \\
Left rotation & LL or LR & $k > s$ & $(r,s)$ \\
Right rotation & LL or LR & $k > s$ & $(r,s)$ \\
 & & & \\
\hline
\end{tabular}
\vspace{.1in}
\caption{Change in $G(c)=(r,s)$ when a single rotation is performed at level $k$ along the right arm of the tree, and the spinal ancestor of $c$ is on the right arm of the tree.  \label{tab:right-G(c)}}
\end{table}

To give a lower bound on the restricted right-arm distance between the trees $T_1$ and $T_2$ which form the tree pair diagram for $w$, we consider the sibling pairings involving leaf $n-2$.  Leaf $n-2$ is paired with leaf $n-3$ in $T_1$ and paired with leaf $n-1$ in $T_2$.  In the following lemmas we bound the minimal number of rotations necessary to split these sibling pairs.  Combined, these estimates yield the desired lower bound.  The main tool is the ordered pair $G(c)$, which allows us to track the position of the exposed caret containing leaf $n-2$ relative to the right arm of the tree.

\begin{lemma}\label{buriedsibs}
Let $w=(T_1,T_2) \in F$ have normal form $$x_{m+2} x_{m+3} \cdots x_{n-2} x_{n-3}^{-1} x_{n-4}^{-1} \cdots x_{m+1}^{-1}$$ where $n > m+4$.  
The tree $T'$ resulting from the application of at most $2n-2m-3$ rotations at locations at levels $0$ to $m$ along the right side of the tree to $T_1$ will contain the sibling pair $[n-3,n-2]$.
\end{lemma}

\begin{proof}
We note that by Lemma \ref{affectedsibs}, sibling pairs can be destroyed by a single rotation only when they are connected by one left edge to a node on the right side of the tree at level $m$ or less, or are the leaves of the rightmost caret in the tree.  The  exposed caret $c$ in $T_1$ with leaves $n-3$ and $n-2$ cannot be moved to be the rightmost caret of the tree, as all rotations preserve the natural infix order on the carets.

Thus, until the caret $c$ is connected by a single left edge to the right side of the tree at level at most $m$, and the correct rotation is performed to separate them, leaves $n-3$ and $n-2$ will remain a sibling pair.  
We use the ordered pair $G(c)$ to monitor the position of $c$ relative to the right arm of the tree while performing a series of rotations.  The leaves $n-3$ and $n-2$ will remain sibling pairs until $G(c)=(1,l)$ for some $l \leq m$, when a single rotation can be performed to separate these leaves.

We consider the sequence of trees $S_0 = T_1$, $S_1,S_2, \ldots ,S_t$ resulting from performing a series of $t$ rotations corresponding to a sequence of $t$ generators $g_1g_2 \cdots g_t$.  Each $S_i$ is the result of
applying $g_i$ to $S_{i-1}$.  We trace the images of the caret $c$ through this 
sequence and denote its image in $S_i$ by $c_i$.  While the exposed leaves of each $c_i$ have the same leaf numbers in $S_i$, the entries in $G(c_i)$ may change as a result of each rotation.  The possible changes in $G(c_i)$ are summarized in Table \ref{tab:right-G(c)}.  

We note that it is possible to move the caret $c$ so that its spinal ancestor is on the left arm of the tree.  Since rotations are not allowed along the left arm of the tree, caret $c$ must be returned to a subtree of a right node before the sibling pair $[n-3,n-2]$ can be split.  This will not happen in any minimal length transformation.

We know that initially, $G(c) = (n-m-3,m+1)$, and the sibling pair $[n-3,n-2]$ is not destroyed until after $G(c_i)=(1,l)$, for some appropriate $l \leq m$. The path $\alpha_c$ has labels $LRRR \cdots R$.  These labels remain unchanged as rotations are performed along the right arm of the tree.  As the length of the path $\alpha_{c_i}$ is decreased, labels are removed sequentially from the beginning of this list, but the remaining labels are never changed by rotations along the right arm of the tree.

We see from Table \ref{tab:right-G(c)} that the rotations which reduce the first coordinate of $G(c_i)$ fall into two types.
\begin{enumerate}
\item Rotations which decrease the first coordinate and increase the second.
\item Rotations which decrease the first coordinate and leave the second unchanged.  These can only happen when the initial two labels of $\alpha_{c_i}$ are $LL$.
\end{enumerate}

Along the initial path $\alpha_c$, there are no adjacent left labels.  To create such a pair of labels, in order to perform a reduction of the first coordinate of $G(c_i)$ but leave the second coordinate unchanged, requires the creation of at least one additional caret with label $L$.  While this is easily accomplished, its creation increases the first coordinate of $G(c_i)$.  The resulting rotation then decreases this coordinate with no net change in $G(c_i)$.  So we see that there will never be any rotations of this second type in a minimal sequence of rotations that splits the sibling pair $[n-3,n-2]$ in $T_1$.

To reduce the first coordinate of $G(c) = (n-m-3,m+1)$ to $1$, we will need at least $n-m-2$ rotations, all of the first type listed above.  Each reduction will increase the second coordinate of $G(c)$ by one.  We will need at least $n-m-2$ additional rotations to reduce the second coordinate back to its starting value, without changing the first coordinate.  We then must perform at least one additional rotation to decrease the second coordinate to $m$ before the sibling pair in question can be split.  This gives a minimum of $2n-2m-3$ rotations before the sibling pair $[n-1,n]$ can be destroyed. 
\end{proof}

We make an analogous argument in the lemma below to bound the minimal number of rotations necessary to split the sibling pair $[n-2,n-1]$ in the tree $T_2$.

\begin{lemma}\label{buriedsibs2}
Let $w=(T_1,T_2) \in F$ have normal form $$x_{m+2} x_{m+3} \cdots x_{n-2} x_{n-3}^{-1} x_{n-2}^{-1} \cdots x_{m+1}^{-1}$$ where $n > m+4$.  
The tree $T'$ resulting from the  application of at most $2n-2m-2$ rotations at locations at levels $0$ to $m$ along the right side of the tree to $T_2$ will contain the sibling pair $[n-2,n-1]$.
\end{lemma}

\begin{proof}
We note that in this case, when $c$ is the caret with exposed leaves numbered $n$ and $n+1$, we
have $G(c)=(n-m-3,m+2)$ and to reduce $G(c)$ to $(0,l)$ with $l \leq m$ will take at least
$(n-m-2)+(n-m-2)+2=2n-2m-2$ rotations by the same analysis as in Lemma \ref{buriedsibs}.
\end{proof}

We combine these lemmas to prove Theorem \ref{thm:rra-sharpfull}.

{\it Proof of Theorem \ref{thm:rra-sharpfull}.}
We consider the reduced tree pair diagram $(T_1, T_2)$ corresponding to the element 
$w= x_{m+2} x_{m+3} \cdots x_{n-2} x_{n-3}^{-1} x_{n-2}^{-1} \cdots x_{m+1}^{-1} \in F$, as in the
lemmas above.
If the restricted right-arm rotation distance between these two trees is $d$, then
we consider the sequence of trees $S_0 = T_1$, $S_1,S_2, \ldots ,S_d=T_2$ resulting from performing that series of $d$ rotations to $T_1$ to get $T_2$.
Lemma \ref{buriedsibs} shows that any application of $2n-2m-3$ allowed rotations
to $T_1$ will still result in a  tree with leaves $n-3$ and $n-2$ still paired, so in trees
$S_i$ with $0 \leq i \leq 2n-2m-3$ leaf $n-2$ must be paired with $n-3$.
We consider the tail end of that sequence, and find that Lemma \ref{buriedsibs2}
shows that in trees $S_i$
with $ d-(2n-2m-2) \leq i \leq d$ leaf $n-2$ must be paired with $n-1$.
Since it will take at least one additional rotation to change the
pairing of leaf $n-2$ from $n-3$ to $n-1$,  the restricted right arm rotation distance
between the two trees is at least $4n-4m-4$.  
\qed

Theorem \ref{thm:rra-sharpfull} gives a family of  pairs of trees with $n$ nodes satisfying a lower bound on restricted right-arm rotation distance 
with respect to a generating set ${\mathcal S'}$ which includes all generators from $x_0$ to $x_m$.  Restricting the generating set to a subset ${\mathcal S}$ of ${\mathcal S'}$ which includes $x_0$ can only increase the restricted right-arm rotation distance between two trees, or cause it to be undefined.  In the case of the words used in the proof of Theorem \ref{thm:rra-sharpfull}, it follows from Lemma \ref{lemma:notreducednotdefined} that the restricted right-arm rotation distance will still be defined when the generating set is further restricted.  This follows because the smallest index in the normal form of the words used exceeds the highest level along the right arm of the tree where rotation is allowed.  Thus we have proven Theorem \ref {thm:rra-sharp} as well.

\section{Bounding right-arm rotation distance}

The original arguments of Culik and Wood \cite{cw} which give a
bound on ordinary rotation distance apply to right-arm rotation
distance as well.   Their argument is that any binary tree $T$
with $n$ nodes can be transformed to or from the all-right tree
with $n$ nodes  by no more than $n-1$
rotations, all of which can be chosen to lie on the right arm of the tree.
Thus, the right-arm rotation distance between two trees $T_1$ and
$T_2$ each with $n$ nodes is no more than $2n-2$, as we can
transform $T_1$ to the all-right tree and from there transform it to
$T_2$.  While this bound is not optimal for the original rotation
distance, we show that it is optimal for right-arm rotation distance.

\begin{theorem}
For each $n \geq 3$, there are rooted binary trees $T_1$ and $T_2$
each with $n$ nodes so that the right-arm rotation distance
between them satisfies $d_{RA}(T_1,T_2) = 2n -2$.
\end{theorem}

\begin{figure}\includegraphics[width=3in]{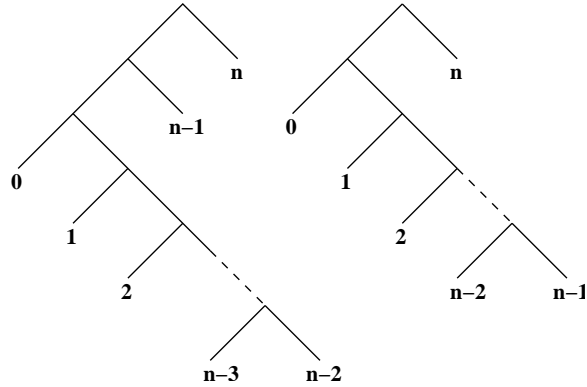}\\
\caption{The tree pair diagram for words of the form $x_0 x_1 x_2
x_3 \cdots x_{n-2} x_{n-3}^{-1} x_{n-4}^{-1} \cdots x_1^{-1}
x_0^{-2}$, with  $n$ nodes and length $2n-2$ with respect to the
infinite generating set.}
\label{longrapic}
\end{figure}

\begin{proof}
To prove this we consider the elements of $F$ with normal form 
$$x_0 x_1 x_2 x_3 \cdots
x_{n-2} x_{n-3}^{-1} x_{n-4}^{-1} \cdots x_1^{-1} x_0^{-2}$$
pictured in Figure \ref{longrapic}, which have $n$ nodes and
have word length $2n-2$
with respect to the infinite generating set for $F$.  It follows from Proposition \ref{prop:infgenset} that this is also the right-arm rotation distance between the two trees.
\end{proof}

\section{Left-arm and spinal rotation distances}

We now consider rotation distances which include rotations at
nodes along the left side of the tree,  instead of or in addition to,
nodes along the right side of the tree.  It is clear by symmetry that {\em
restricted left-arm rotation distance}, which allows rotations only at a finite collection of nodes on the left side of the tree and the root node, will satisfy
the same  bounds as restricted right-arm rotation distance.
Similarly, {\em left-arm rotation distance}, which allows rotations at any node along the left arm of the tree, will satisfy the sharp upper bound of $2n-2$ on trees with $n$ nodes, for $n \geq 3$.

Finally, we consider a rotation distance which allows rotations at the root
node, a finite nonempty collection of nodes on the right side of
the tree, and at a finite nonempty collection of nodes on the left
side of the tree.  Since all nodes where rotations are permitted
lie on the spine of the tree, we call such a rotation distance
a {\em restricted spinal rotation distance}.  
In terms of Thompson's group $F$, left rotation at level
$n$ on the left arm of the tree can be expressed as $y_n = x_0^n x_1 x_0^{-n-1}$.
So there is an {\em extended infinite generating set} for $F$ consisting of all
$x_n$ and $y_n$ which corresponds to allowing rotations at
any location on the spine.  Here, we consider finite subsets of this
enlarged generating set.  Again, if we do not include $x_0$ in the
generating set we consider, we generate either a subgroup isomorphic
to $F$ or its direct square, so we restrict to the case where $x_0$ is included
in the generating set.

\begin{definition}
Let ${\mathcal S}=\{x_0,x_{i_1}, \ldots ,x_{i_L}, y_{j_1}, \ldots, y_{j_l} \}$ with $i_1 < i_2 \cdots < i_L$ and $j_1 < j_2 \cdots < j_l$ be a finite subset of the extended infinite generating set for $F$
and $T_1$ and $T_2$ be trees with the same number of leaves.
We define  $d_{RS}^{\mathcal S} (T_1,T_2)$,  the {\em  restricted spinal rotation distance
with respect to  ${\mathcal S}$}, as the minimal number of rotations required to transform $T_1$ to
$T_2$, where the rotations are only allowed at levels $0, i_1, \ldots, i_{L-1}$ and
$i_L$  along
the right side of the tree and at levels $j_1, \ldots, j_l$ on the left side of the tree.

\end{definition}

Again, though allowing rotations at finitely many
locations on both the right and left arms  of the tree may reduce the rotation distance between some pairs of trees, we prove that the multiplicative
constant of 4 in the upper bound cannot be decreased. 

We first show that spinal rotation distance satisfies the same upper bound as restricted right-arm rotation distance.

\begin{theorem}\label{thm:spinal-upperbound}
Let ${\mathcal S}=\{x_0,x_{i_1}, \ldots ,x_{i_L}, y_{j_1}, \ldots, y_{j_l} \}$ with $i_1 < i_2 \cdots < i_L$ and $j_1 < j_2 \cdots < j_l$, where $x_i$ is a generator of $F$ and $y_n = x_0^{n} x_1 x_0^{-n-1}$, and let $d_{RS}^{\mathcal S}$ be the corresponding spinal rotation distance.  Let $T_1$ and $T_2$ be binary trees, each with $n$ carets with $n \geq 3$,
 for which $d_{RS}^{\mathcal S}(T_1,T_2)$ is defined.  Then
$$d_{RS}^{\mathcal S}(T_1,T_2) \leq 4n-8.$$
\end{theorem}

\begin{proof}
Let ${\mathcal S'} = \{x_0,x_{i_1}, \ldots ,x_{i_L} \}$, which is also a generating set for $F$ in which each generator corresponds to a rotation along the right arm of the tree.  Since adding additional elements to a generating set can only decrease the corresponding restricted rotation distance, we see immediately that $d_{RS}^{\mathcal S}(T_1,T_2) \leq d_{RS}^{\mathcal S'}(T_1,T_2)$. Since $d_{RS}^{\mathcal S'}(T_1,T_2) = d_{RRA}^{\mathcal S'}(T_1,T_2)$, and Proposition \ref{prop:rrabounds} proves that $d_{RRA}^{\mathcal S'}(T_1,T_2) \leq 4n-8$, the theorem follows.
\end{proof}

Now we show that the multiplicative coefficient of 4 is optimal in the same sense as with restricted right-arm rotation distance.  
For these examples, to avoid possible repeated excessive reductions in $G(c)$, we take words which have an exposed caret connected to the right-hand side of the tree with a path which alternates between branching right and left.

\begin{theorem}\label{thm:rra-sharpboth}
Let ${\mathcal S}=\{x_0,x_{i_1}, \ldots ,x_{i_L}, y_{j_1}, \ldots, y_{j_l} \}$ with $i_1 < i_2 \cdots < i_L$ and $j_1 < j_2 \cdots < j_l$, where $x_i$ is a generator of $F$ and $y_n = x_0^{n} x_1 x_0^{-n-1}$.  Then there exist trees 
$T_1$ and $T_2$ with $n$  nodes where $n > \max \{i_L,j_l\}$ for which $d_{RS}^{\mathcal S}(T_1,T_2)$ is defined that satisfy 
 $$ d_{RS}^{\mathcal S}(T_1,T_2) \geq   4n-4 \max \{i_L,j_l\} - 12   .$$
\end{theorem}

\begin{figure}\includegraphics[width=3.5in]{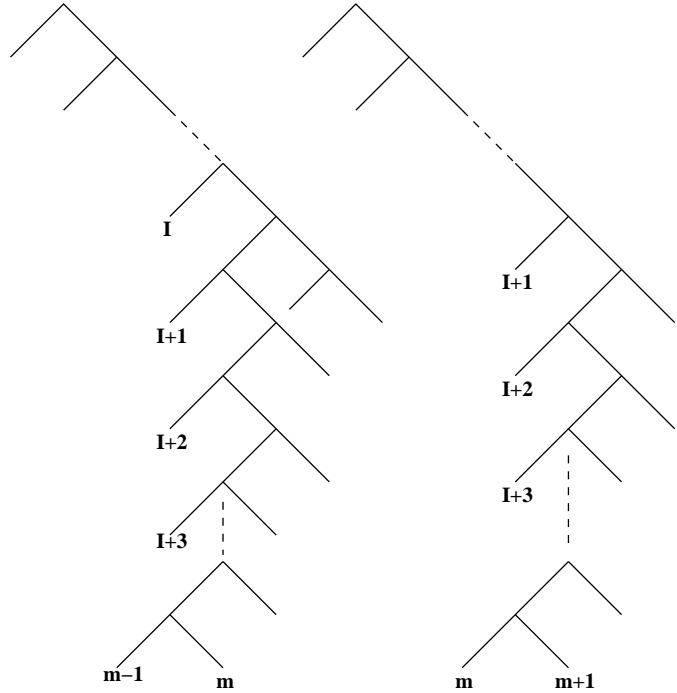}\\
\caption{An example from the family of elements of the form $w = (T_1,T_2) =
x_{I+2} x_{I+3}^2 x_{I+4}^2 \cdots x_{m}^2 x_{m-1}^{-2} x_{m-2}^{-2}  \cdots x_{I+2}^{-2} x_{I+1}^{-1}
 $ which are used to show that the multiplicative constant of $4$ in the upper bound on restricted spinal rotation distance is optimal.}
\label{fig:spinalword}
\end{figure}

We again introduce a particular family of elements $w \in F$, represented by reduced tree pair diagrams $(T_1,T_2)$, which requires this lower bound on the restricted spinal
rotation distance between $T_1$ and $T_2$.  As in the proof of Theorem \ref{thm:rra-sharpfull}, we use the ordered pair $G(c)$ as the main tool of the proof.  Since rotations are now permitted along the left arm of the tree, we must note the changes in $G(c)$ caused by a rotation along the left arm of the tree.  As before, $c$ is an exposed caret, $\alpha_c$ is the minimal path from the node of $c$ to the spinal ancestor of $c$, each non-spinal node along $\alpha_c$ is given a label of $R$ or $L$, and $G(c) = (r,s)$.  If the spinal ancestor of $c$ lies on the left arm of the tree, then the initial label along $\alpha_c$ must be $R$.

The following table summarizes the changes in $G(c)$ when a single rotation is performed along the left arm of the tree.  If the spinal ancestor of $c$ is on the right arm of the tree, and a rotation is performed on the left arm of the tree, at level at least one, then $G(c)$ remains unchanged.

\begin{table}[h]\label{G(c)summary-left}
\begin{tabular}{|c|c|c|c|}
\hline
Direction & Initial labels  & Relative position & Change in \\
of Rotation  & of $\alpha_c$ &  of $k$ and $s$ & $G(c) = (r,s)$  \\
\hline 
 & & & \\
Left rotation & RL or RR & $k < s-1$ & $(r,s+1)$ \\
Right rotation & RL or RR & $k < s-1$ & $(r,s-1)$ \\ 
 & & & \\
\hline
 & & & \\
Left rotation & RL or RR & $k = s-1$ & $(r,s+1)$ \\
Right rotation & RL or RR & $k = s-1$ & $(r+1,s-1)$ \\
 & & & \\
\hline
 & & & \\
Left rotation & RL & $k = s$ & $(r-1,s+1)$ \\
Left rotation & RR & $k = s$ & $(r-1,s)$ \\
Right rotation & RL or RR & $k = s$ & $(r+1,s)$ \\
 & & & \\
\hline
 & & & \\
Left rotation & RL or RR & $k > s$ & $(r,s)$ \\
Right rotation & RL or RR & $k > s$ & $(r,s)$ \\
 & & & \\
\hline
\end{tabular}
\vspace{.1in}
\caption{Change in $G(c)=(r,s)$ when a single rotation is performed at level $k \geq 1$ along the left arm of the tree, and the spinal ancestor of $c$ is on the left arm of the tree.  \label{tab:left-G(c)}}
\end{table}

Rotation at the root caret, corresponding to multiplication by the generator $x_0^{\pm 1}$ can affect $G(c)$ in one of two ways.
\begin{enumerate}
\item If $G(c) = (r,1)$, then this rotation may change the arm of the tree on which the spinal ancestor of $c$ lies.  If this is the case, then the first label along the path $\alpha_c$ will change from $L$ to $R$ or vice versa and $G(c)$ will remain $(r,1)$.
\item If $G(c) \neq (r,1)$ then the $r$ coordinate remains unchanged, and the $s$ coordinate is changed by $\pm 1$, depending on the direction of the rotation.
\end{enumerate} 
Other rotations which do not add carets to the path $\alpha_c$ will not change the labels along $\alpha_c$.  Rotations which decrease the length of $\alpha_c$ can only remove the initial label along the path.

The proof of Theorem \ref{thm:rra-sharpboth} uses methods analogous to the proofs of Lemmas \ref{buriedsibs} and \ref{buriedsibs2} and Theorems \ref{thm:rra-sharp} and \ref{thm:rra-sharpfull}.  

{\it Proof of Theorem \ref{thm:rra-sharpboth}.} First, we can assume without loss of generality
that $i_L$ is larger than $j_l$, since if not, we can interchange the left and right sets of generators
by taking reflections of the trees considered.
Given the set ${\mathcal S}$, we define the level furthest from the root on the right side at which a rotation can take place as $I=i_L$ for convenience and
consider the elements with $m>I$ of the form 
$$w=(T_1,T_2) = x_{I+2} x_{I+3}^2 x_{I+4}^2 \cdots x_{m}^2 x_{m-1}^{-2} x_{m-2}^{-2}  \cdots x_{I+2}^{-2} x_{I+1}^{-1}.$$

The reduced tree pair diagram $(T_1,T_2)$ for $w$ is pictured in
Figure \ref{fig:spinalword}. 
 Each $T_i$ contains $n = 2m-I$ carets, with a deeply buried exposed caret in each tree
 connected to the right-hand side of the tree via a zigzag path.
 In $T_1$, we label the caret with the sibling pair $[m-1,m]$ as $c$
 and in $T_2$, we label caret with the sibling pair $[m,m+1]$ as $d$.
  We see that the lengths of $\alpha_c$ and $\alpha_d$ are both $2m-2I-3$, and have labels $LRLRLRL \cdots RL$.  We begin with $G(c) = (2 m-2I-3,I+1)$ and $G(d)=(2m-2I-3,I+2)$ and to split each sibling pair, we need to reduce both $G(c)$ and $G(d)$ to $(1,l)$ for some appropriate $l$.
As in the right-arm case described above, when counting the needed rotations, we count from
the beginning of the transformation to determine the number of rotations needed to separate the
sibling pair in caret $c$ and from the end of the transformation to determine the number of rotations to
separate the pair in caret $d$.

As in the proof of Theorem \ref{thm:rra-sharpfull}, there are two ways that a single rotation can decrease the $r$ coordinate of $G(c)$, which we extract from Tables \ref{tab:right-G(c)} and \ref{tab:left-G(c)}.
\begin{enumerate}
\item Rotations which decrease the first coordinate and increase the second.
\item Rotations which decrease the first coordinate and leave the second unchanged.  These can only happen when the spinal ancestor of $c$ is on the right arm of the tree and initial two labels of $\alpha_{c_i}$ are $LL$ or when the spinal ancestor of $c$ is on the left arm of the tree and initial two labels of $\alpha_{c_i}$ are $RR$.  We will call these {\em bonus} rotations.
\end{enumerate}

We note that a single rotation which does not increase the $r$ coordinate either leaves the labels along $\alpha_c$ unchanged, or removes the initial label.  There is no way for a single rotation to change a label in the middle of the path.  Rotation at the root may change the initial label from $R$ to $L$ or vice versa, or leave all labels unchanged.

We first enlarge our generating set to ${\mathcal S'} = \{x_0,x_1,x_2, \cdots x_{i_L},y_1,y_2, \cdots ,y_{j_l}\}$, so that rotations are permitted at all nodes at levels zero through $i_L$ along the right side of the tree, and at all levels one through $j_l$ along the left side of the tree.  As with Theorems \ref{thm:rra-sharp} and \ref{thm:rra-sharpfull}, if we can produce trees $T_1$ and $T_2$ which require the desired lower bound on the spinal rotation distance $d_{RS}^{\mathcal S'}(T_1,T_2)$, then the same bound holds with respect to ${\mathcal S}$, since the removal of elements from the generating set can only increase the spinal rotation distance or cause it to be undefined.  Since the elements in our example do not have any generators in their normal forms of index less than $i_L$, Lemma \ref{lemma:notreducednotdefined} guarantees that the relevant restricted right-arm rotation distance using just the rotations on the right side of the tree is defined.  Thus the restricted spinal rotation distance using the larger set of rotations corresponding to the entire generating set ${\mathcal S}$ will also be defined.

We now give a lower bound on the minimal number of rotations necessary to separate the sibling pair $[m-1,m]$ in $T_1$, which we can equivalently view as multiplication by a minimal sequence of generators.  We recall that we must reduce $G(c)$ from $(2 m-2I-3,I+1)$, where $I = i_L$, to $(1,l)$ for some level $l$ at which rotation is allowed on the correct arm of the tree.  Every rotation, whether at a node on the right or left arm of the tree, which decreases the $r$ coordinate of $G(c)$, with the exception of  bonus rotations, also increases the $s$ coordinate.  The $s$ coordinate of $G(c)$
begins larger than $i_L$ (and  thus $j_l$ as well), and so each non-bonus rotation
which reduces the $r$ coordinate will require an additional rotation to counteract
the corresponding increase in the $s$ coordinate.

We now consider the role of bonus rotations that can be used in a minimal length
sequence of rotations.  We note that the labels for $\alpha_c$ begin with $LRLRLRL\cdots$.
There is no natural occurrence in the sequence of
labels of $RR$ or $LL$.  In order to create an $RR$ or $LL$, we can use an application
of  at least one $x_0^{\pm1}$ which will change the initial label but will not change $G(c)$.
Thus each such bonus rotation would need to be
accompanied by an  $x_0^{\pm1}$ which does not change $G(c)$.  

Let $\beta$ be a minimal string of rotations which reduces $G(c) = (2m-2I-3,I+1)$ to $(1,l)$.  We divide the rotations in $\beta$ into three groups:
\begin{enumerate}
\item $p_1$ non-bonus rotations which change the coordinates of $G(c)$,
\item $p_2$ bonus rotations which change the coordinates of $G(c)$, and 
\item $p_3$ rotations which do not change $G(c)$.
\end{enumerate}

The argument above showing that there is at least one
$x_0^{\pm1}$ accompanying each bonus rotation 
shows that $p_3 \geq p_2$.  Totalling the effects of these $p$ rotations, we find that to reduce $G(c) = (2m-2I-3,I+1)$ to $G(c) = (1,l)$ with $l \leq I$
will require at least $2(2m-2I-4)+1= 2n-2I-7$ total rotations, since $n = 2m-I$.  Thus, as in Lemma \ref{buriedsibs}, any sequence
of at most $2n-2I-7$ rotations applied to tree $T_1$ will have leaves $m-1$ and $m$ as sibling
pairs.

Similarly, working backward and considering tree $T_2$, we need to change
the $G(d)=(2m-2I-3,I+2)$ to $(1,l)$ with $l \leq I$, to be able to affect the sibling pair
$[m,m+1]$.  A similar calculation shows that
 any sequence
of at most $2n-2I-6$ rotations applied to the tree $T_2$ will have leaves $m$ and $m+1$ as sibling
pairs.

So we see that 
in the sequence of trees $S_0 = T_1, S_1,S_2, \ldots ,S_d=T_2$ exhibiting
the transformation of $T_1$ into $T_2$ via rotation, the sibling pair $[m-1,m]$ must be present in the first $2n-2I-7$
trees and the sibling pair $[m,m+1]$ must be present in the last $2n-2I-6$ trees.  
Including the rotation to change the sibling pairing of leaf $m$, we see
that the spinal rotation distance must be at least $4n-4I-12$.

Note that the number of nodes in these trees $n$ is constructed to be $2m-I$, so depending
upon whether or not $I$ is even or odd, the number of nodes in the trees constructed using these examples for increasing $m$ will either always be even or always be odd.  To obtain examples for all parity $n$ larger than $I$, we
can repeat the argument above on examples with one additional caret, with normal
forms 
$$w'=(T_1,T_2) = x_{I+2} x_{I+3}^2 x_{I+4}^2 \cdots x_{m}^2 x_{m+1}  x_m^{-1}  x_{m-1}^{-2} x_{m-2}^{-2}  \cdots x_{I+2}^{-2} x_{I+1}^{-1}$$
and find that the sibling pair $[m,m+1]$ must be present in the first $2n-2I-5$ steps
and that the sibling pair $[m+1,m+2]$ must be present in the last $2n-2I-6$ steps,
so in those cases the distance between the two trees must be at least $4n-4I-10$.
Thus, we have that the bound holds for all $n$ larger than $I$ and thus the therorem.
\qed

\bibliographystyle{plain}

\end{document}